\documentclass[aos,MSNbibl,nameyear,secthm,seceqn,dvips]{arximspdf}

\doi{10.1214/11-AOS875}
\volume{39}
\issue{2}
\pubyear{2011}
\firstpage{1266}
\lastpage{1281}

\makeatletter

\def\bptnote#1{}

\newcommand{\eqref}[1]{(\ref{#1})}

\newtheorem{lem}[thm]{Lemma}

\newproclaim{remark}[thm]{Remark}

\makeatother

\begin{document}
\begin{frontmatter}

\title{A note on the de la Garza phenomenon for~locally~optimal designs}
\runtitle{de la Garza phenomenon for locally optimal designs}

\begin{aug}
\author[A]{\fnms{Holger} \snm{Dette}\corref{}\ead[label=e1]{holger.dette@rub.de}\thanksref{t1}}
and
\author[B]{\fnms{Viatcheslav B.} \snm{Melas}\ead[label=e2]{vbmelas@post.ru}\thanksref{t2}}
\thankstext{t1}{Supported in
part by the Collaborative Research Center ``Statistical modeling of
nonlinear dynamic processes'' (SFB 823) of the German Research Foundation (DFG).}
\thankstext{t2}{Supported in part by Russian Foundation for Basic
Research (RFBR), Project \mbox{09-01-00508}.}
\runauthor{H. Dette and V. B. Melas}
\affiliation{Ruhr-Universit\"{a}t Bochum and St. Petersburg State University}
\address[A]{Fakult\"{a}t f\"{u}r Mathematik\\
Ruhr-Universit\"{a}t Bochum \\
44780 Bochum\\
Germany \\
\printead{e1}}

\address[B]{Department of Mathematics\\
St. Petersburg State University\\
St. Petersburg\\
Russia\\
\printead{e2}}
\end{aug}

\received{\smonth{7} \syear{2010}}
\revised{\smonth{1} \syear{2011}}

%
\begin{abstract}
The celebrated de la Garza phenomenon states that for a polynomial
regression model of degree $p-1$ any optimal design can be based on at most
$p$ design points. In a remarkable paper, Yang [\textit{Ann. Statist.}
\textbf{38} (\citeyear{yang2010}) 2499--2524]
showed that this phenomenon exists in many locally optimal design
problems for
nonlinear mo\-dels. In the present note, we present a different view
point on these findings using results about moment theory and Chebyshev
systems. In particular, we show that this phenomenon occurs in an even
larger class of models than considered so far.
\end{abstract}

\setattribute{keyword}{AMS}{\textit{MSC2010} subject classification.}
%

\end{frontmatter}

\section{Introduction}\label{sec1}

Nonlinear regression models are widely used for modeling dependencies
between response and explanatory variables [see \citet{sebwil1989} or
\citet{ratkowsky1990}]. It is well known that an appropriate choice of
an experimental design can improve the quality of statistical analysis
substantially, and therefore the problem of constructing optimal
designs for nonlinear regression models has found considerable
attention in
the literature. Most authors concentrate on locally optimal designs
which assume that a guess for the unknown parameters of the model is
available [see \citet{chernoff1953}, \citet{fortorwu1992}, \citet
{hestusun1996}, \citet{fangheda2008}]. These designs are usually used as
benchmarks for commonly used designs. Additionally, they serve as a
basis for constructing optimal designs with respect to more sophisticated
optimality criteria which address for a less precise knowledge about
the unknown parameters [see \citet{pronwalt1985} or \citet{chaver1995},
\citet{dette1997}, \citet{muepaz1998}].
It is a well-known fact that the numerical or
analytical calculation of optimal designs simplifies substantially if
it is known that the optimal design is saturated, which means that the
number of different experimental conditions coincides with the number
of parameters in the model [see, e.g., \citet{hestusun1996},
\citet{dettwong1996}, \citet{imhostud2001}, \citet{imhof2001}, \citet
{melas2006}, \citet{fangheda2008}
among many others].

So, the ideal situation appears if the optimal design is in the
sub-class of all saturated designs. In a celebrated paper, \citet{garza1954}
proved that for a $(p-1)$th-degree
polynomial regression model, any optimal design can be based on at most
$p$ points. \citet{khumuksingho2006} considered a nonlinear regression
model and introduced the terminology of the de la Garza phenomenon,
which means that for any design there exists a saturated design, such that
the information matrix of the saturated design is not inferior to that
of the given design under the Loewner ordering. In a remarkable paper,
\citet{yang2010} derived sufficient conditions on the nonlinear
regression model for the occurrence of the de la Garza phenomenon and
demonstrated that this situation appears in a broad class of nonlinear
regression models. These results generalize recent findings of
\citet{yangstuf2009} for nonlinear models with two parameters.

However, some care is necessary if these results are applied as
indicated in the following simple example of homoscedastic linear
regression on
the interval $[0,1]$.
Here the information matrix of the design which advises the
experimenter to take all $n$ observations at the point 0 is given by\vspace{-2pt}
\[
X^T_1 X_1
=
\pmatrix{ n & 0 \cr 0 & 0
}
\]
while any other design (using the experimental conditions
$x_1,\dots,x_n$) yields an information matrix\vspace{-3pt}
\[
X^T_2 X_2 =
\pmatrix{ n &\displaystyle  \sum^n_{i=1}x_i \cr \displaystyle \sum^n_{i=1}x_i &
\displaystyle \sum^n_{i=1}x_i^2
}.
\]
It is easy to see that the matrix $X^T_2 X_2 - X^T_1 X_1$ is
indefinite (i.e., it has positive
and negative eigenvalues) whenever one of the $x_i$ is positive.
Consequently, the design corresponding to $X^T_1 X_1$ cannot be
improved. On
the other hand, it is also easy to see that for any $k \in\{1,\ldots,
\lfloor n/2 \rfloor-1 \} $ the information matrix of the design, which
takes observations at $x_1=\cdots= x_{n-2k}=0$ and at $x_{n-2k+1}=
\cdots= x_{n}=1/2$ can be improved (with respect to the Loewner ordering)
by the information matrix corresponding to the design $x_1=\cdots
=x_{n-k}=0$ and \mbox{$x_{n-k+1}= \cdots=x_{n}=1$}. Thus, there exist designs where
a ``real'' improvement is possible, while other designs cannot be
improved. Note that the results in \citet{yang2010} do not provide a
classification of the two types of designs.

It is the purpose of the present paper to present a more detailed view
point on these problems, which clarifies this---on a first glance---contradiction. In contrast to the method used by \citet{yang2010},
which is mainly algebraic, our approach is analytic and based on the theory
of Chebyshev systems and moment spaces [see \citet{karstu1966}]. In
particular, we will demonstrate that the de la Garza phenomenon appears in
any nonlinear regression model, where the functions in the Fisher
information matrix form a Chebyshev system. Additionally, we will solve the
problem described in the previous paragraph and we will
identify the sufficient conditions stated in
\citet{yang2010} as a special case of an extended Chebyshev system.
Therefore, our results generalize the recent findings of \citet{yang2010} in
a nontrivial way and, additionally, provide---in our opinion---a more
transparent and more complete explanation of the de la Garza phenomenon
for optimal designs in nonlinear regression models.

The remaining part of this paper is organized as follows. Section \ref{sec2}
provides a brief introduction in the problem, while Section \ref{sec3} contains our
main results. Finally, the new results are illustrated in a rational
regression model, where the currently available methodology cannot be used
to establish the de la Garza phenomenon.

\section{Locally optimal designs}\label{sec2}

Consider the common nonlinear regression model\vspace{-2pt}
%
\begin{equation}\label{1.1}
Y=\eta(x,\theta)+\varepsilon,
\end{equation}
where $\theta\in\Theta\subset\mathbb{R}^p$ is the vector of unknown
parameters, and different observations are assumed to be independent. The
errors are normally distributed with mean $0$ and variance $\sigma^2$.
The variable $x$ denotes the explanatory variable, which varies in the
design space $[A,B] \subset\mathbb{R}$. We assume that $\eta$ is a
continuous and real valued function of both arguments $(x,\theta) \in[A,B]
\times\Theta$ and differentiable with respect to the variable $\theta
$. A design is defined as a~probability measure $\xi$ on the interval
$[A,B] $ with finite support [see \citet{kiefer1974}]. If the design
$\xi$ has masses $w_i$ at the points $x_i$ $(i = 1, \dots, k)$ and $n$
observations can be made by the experimenter, this means that the
quantities $w_i n$ are rounded to integers, say $n_i$, satisfying
$\sum^k_{i=1} n_i =n$, and the experimenter takes $n_i$ observations
at each location $x_i$ $(i=1, \dots, k)$.
The information matrix of an approximate design $\xi$ is defined by\vspace{-2pt}
%
\begin{equation} \label{1.1a}
M(\xi, \theta) = \int^B_A \biggl( \frac
{\partial}{\partial\theta} \eta(x,\theta)\biggr) \biggl( \frac{\partial}{\partial
\theta} \eta(x,\theta) \biggr)^T \,d \xi(x),
\end{equation}
and it is well known [see
\citet{jennrich1969}] that under appropriate assumptions of
regularity the covariance matrix of the least squares estimator is
approximately given by $\sigma^2 \ M^{-1}(\xi,\theta)/n$, where $n$ denotes
the total sample size and we assume that the observations are taken
according to the approximate design $\xi$.

An optimal design maximizes an appropriate functional of the
information matrix and numerous criteria have been proposed in the
literature to
discriminate between competing designs [see \citet{silvey1980}, \citet
{pazman1986} or \citet{pukelsheim2006} among others]. Note that in
nonlinear regression models the information matrix (and as a
consequence the corresponding optimal designs) depend on the unknown parameters
and are therefore called locally optimal designs [see \citet
{chernoff1953}]. These designs require an initial guess of the unknown
parameters
in the model and are used as benchmarks for many commonly used designs.

Most of the available optimality criteria satisfy a monotonicity
property with respect to the Loewner ordering, that is
%
\begin{equation}\label{low}
M(\xi_1, \theta) \leq M(\xi_2,\theta) \quad \Longrightarrow\quad \Phi(M(\xi
_1,\theta)) \leq\Phi(M(\xi_2,\theta)),
\end{equation}
where the parameter
$\theta$ is
fixed, $\xi_1, \xi_2$ are two competing designs and $\Phi$ denotes
an information function in the sense of \citet{pukelsheim2006}.
For this reason, it is
of interest to derive a complete class theorem in this general context
which characterizes the class of designs, which cannot be improved with
respect to the Loewner ordering of their information matrices. We call
a design $\xi_1$ admissible if there does not exist a design $\xi_2$,
such that $M(\xi_1,\theta)\neq M(\xi_2,\theta)$ and
%
\begin{equation}\label{adm}
M( \xi_1,\theta) \leq M(\xi_2,\theta).
\end{equation}
As pointed out in
\citet{yang2010} for many nonlinear regression models the information
matrix defined in \eqref{1.1a} has a representation of the form
%
\begin{equation}\label{rep1}
M(\xi,\theta) = P(\theta) C (\xi,\theta) P^T (\theta),
\end{equation}
where $P(\theta)$ is a nonsingular $p \times p$ matrix, which
does not
depend on the design~$\xi$, the matrix $C$ is defined by
%
\begin{equation} \label{cmat}
C(\xi,\theta) = \pmatrix{
\displaystyle \int^B_A \Psi_{11} (x) \,d \xi(x) & \cdots& \displaystyle \int^B_A \Psi_{1p} (x)
\,d \xi(x) \cr
\vdots& \ddots& \vdots\cr
\displaystyle \int^B_A \Psi_{p1} (x) \,d \xi(x) & \cdots& \displaystyle \int^B_A \Psi_{pp} (x)
\,d \xi(x)
}
\end{equation}
and $\Psi_{11},\Psi_{12},\dots,\Psi_{pp}$ are functions defined on
the interval $[A,B]$. Note that these functions usually depend on the
parameter $\theta$, but for the sake of simplicity we do not reflect
this dependence in our notation. Obviously the inequality \eqref{adm} is
satisfied if and only if the inequality
%
\begin{equation}\label{adm1}
C(\xi_1,\theta) \leq C(\xi_2,\theta)
\end{equation}
is satisfied.

\section{Chebyshev systems and complete class theorems}\label{sec3}

In the following discussion, we make extensive use of the property that
a system of functions has the Chebyshev property. Following
\citet{karstu1966},
a set of $k+1$ continuous functions $u_0, \ldots, u_k \dvtx [A,B] \to
\mathbb{R}
$ is called a Chebyshev system
(on the interval $ [A,B]$) if the inequality
%
\begin{equation}
\label{3.11}
\left\vert
\matrix{
u_0(x_0) &u_0(x_1) &\ldots&u_0(x_k) \cr
u_1(x_0) &u_1(x_1) &\ldots&u_1(x_k) \cr
\vdots&\vdots& \ddots&\vdots\cr
u_k(x_0) &u_k(x_1) &\ldots&u_k(x_k)
}
\right\vert >  0
\end{equation}
holds for all $A \le x_0 < x_1 < \cdots< x_k \le B$. Note that if the
determinant in~\eqref{3.11} does not vanish then either the functions
$u_0, u_1, \ldots,u_{k-1},u_k$ or the functions $u_0, u_1, \ldots
,u_{k-1}, -u_k$ form a Chebyshev system. The Chebyshev property has widely
been used to determine explicitly $c$-optimal designs [see \citet
{hestusun1996}, \citet{detmelpepstr2003} or \citet{debrpepi2008} among many
others]. On the other hand, its application to other optimality
criteria has not been studied intensively. In the following discussion,
we will
demonstrate that this property will essentially be the reason for the
occurrence of the de la Garza phenomenon. In particular, we will show
that it is essentially sufficient to obtain a complete class theorem
for the design problems associated with the nonlinear regression model~\eqref{1.1}.

For this purpose, we define the index $I(\xi)$ of a
design $\xi$ on the interval $[A,B]$ as the number of support points, where
the boundary points $A$ and $B$ (if they occur as support points) are
only counted by $1/2$. Recall the definition of the matrix $C$ in
\eqref{cmat} and denote by $\Psi_1,\dots,\Psi_k$ the different
elements among the functions $\{ \Psi_{ij} \mid1 \leq j,j \leq p \}$,
which are
not equal to the constant function. Throughout this paper, we assume
\begin{eqnarray} \label{def}
\Psi_k &=& \Psi_{ll}  \qquad \mbox{for some }   l \in\{1,\dots,p\}
\quad \mbox{and} \nonumber\\[-8pt]\\[-8pt]
\nonumber
\Psi_{ij} &\neq&\Psi_{k}\hspace*{1.26pt} \qquad \mbox{for all } (i,j) \neq(l,l)
\end{eqnarray}
[see \citet{yang2010}].
Additionally, we put $\Psi_0(x)=1$ and assume either that
\begin{eqnarray} \label{cheb+}
&& \{ \Psi_0, \Psi_1, \dots, \Psi_{k-1} \} \quad \mbox{and} \nonumber\\[-8pt]\\[-8pt]
\nonumber
&& \{ \Psi_0, \Psi_1, \dots, \Psi_{k-1}, \Psi_{k} \}
\end{eqnarray}
are
Chebyshev systems
or that
%
\begin{eqnarray} \label{cheb-}
&& \{ \Psi_0, \Psi_1, \dots, \Psi_{k-1} \} \quad \mbox{and} \nonumber\\[-8pt]\\[-8pt]\nonumber
&& \{ \Psi_0, \Psi_1, \dots, \Psi_{k-1}, - \Psi_{k} \}
\end{eqnarray}
are Chebyshev systems then the following result characterizes the class
of admissible designs.

\begin{thm}\label{thm3.1}
\textup{(1)} If the functions $\Psi_0 (x)=1,\Psi_1,\ldots,\Psi
_{k-1},\Psi_k$ satisfy \eqref{def} and \eqref{cheb+}, then for any
design $\xi$
there exists a design $\xi^+$ with at most $\frac{k+2}{2}$ support
points, such that $M(\xi^+, \theta) \geq M(\xi,\theta)$. If the
index of
the design~$\xi$ satisfies
\[
I(\xi) < {k \over2}
\]
then the design $\xi^+$ is uniquely determined in the class of all
designs $\eta$ satis\-fying
%
\begin{equation}\label{unique}
\int^B_A \Psi_i (x) \,d \eta(x) = \int^B_A \Psi_i(x)\,d\xi(x),\qquad
i=0,\dots, k-1,
\end{equation}
and coincides with the design $\xi$. Otherwise [in the case $I(\xi)
\ge{k \over2}$], the following two assertions are valid.
\begin{itemize}[(1a)]
\item[(1a)] If $k$ is odd, then $ \xi^+$ has at most $\frac{k+1}{2}$
support points and $ \xi^+$ can be chosen
such that its support contains the point $B$.\vspace*{1pt}
\item[(1b)] If $k$ is even, then $\xi^+$ has at most $\frac{k}{2}+1$
support points and $ \xi^+$ can be chosen
such that the support of $\xi^+$ contains the points $A$ and $B$.
\end{itemize}
{\smallskipamount=0pt
\begin{longlist}[(2)]
\item[(2)] If the functions $\Psi_0 (x)=1,\Psi_1,\ldots,\Psi
_{k-1},\Psi_k$ satisfy \eqref{def} and \eqref{cheb-}, then for any
design $\xi$
there exists a design $\xi^-$ with at most $\frac{k+2}{2}$ support
points, such that $M(\xi^-, \theta) \geq M(\xi,\theta)$. If the
index of
the design $\xi$ satisfies
\[
I(\xi) < {k \over2}
\]
then the design $\xi^-$ is uniquely determined in the class of all
designs $\eta$ satisfying \eqref{unique} and coincides with the design
$\xi$. Otherwise [in the case $I(\xi) \ge{k \over2}$], the
following two assertions are valid.
\begin{itemize}[(2a)]
\item[(2a)] If $k$ is odd, then $\xi^-$ has at most $\frac{k+1}{2}$
support points and $\xi^-$ can be chosen
such that its support contains the point $A$.
\item[(2b)] If $k$ is even, then $\xi^-$ has at most $\frac{k}{2}$
support points.
\end{itemize}
\end{longlist}}
\end{thm}

\begin{pf}
We only present a proof of the first part (1) of the
theorem, the second part follows by similar arguments. For $i=0,\dots,k$
let
\[
d_i(\xi) = \int^B_A \Psi_i(x)\,d\xi(x)
\]
denote the $i$th ``moment'' and define
\[
\vec{d}_k (\xi)=(d_0(\xi),\dots, d_k (\xi) )^T
\]
as the vector of all ``moments'' up to the order $k$. Consider two
designs $\xi_1$ and $\xi_2$ with
\[
\vec{d}_{k-1} (\xi_1) = \vec{d}_{k-1} (\xi_2)\quad  \mbox{and}\quad  d_k(\xi
_1) \leq d_k(\xi_2),
\]
then for any vector $z=(z_1,\dots,z_p)^T \in\mathbb{R}^p$ we have
for some $l\in\{1,\ldots,p\}$
\[
z^T \bigl(C(\xi_2,\theta)-C(\xi_1,\theta)\bigr) z\geq z^2_l \bigl(d_k(\xi
_2)-d_k(\xi_1)\bigr) \geq0,
\]
which means that
\[
C(\xi_2,\theta) \geq C (\xi_1,\theta).
\]
Now let for a fixed vector of ``moments'' $\vec{d}_{k-1}(\xi)$
\[
d^+_k = \sup\{ d_k(\eta) \mid\eta\mbox{ design on } [A,B] \mbox
{ with } \vec{d}_{k-1} (\eta) = \vec{d}_{k-1} (\xi) \}
\]
denote the maximum of the $k$th ``moment'' over the set of all designs
with fixed ``moments'' up to the order $k-1$. Due to the compactness of
the design space and the continuity of the functions $\Psi_0,\dots
,\Psi_k$, there exists a design~$\xi^+$ such that
%
\begin{eqnarray}
\label{b1} d_j( \xi^+) &=& d_j(\xi) ;\qquad  j=0,\dots,k-1, \\
\label{b2} d_k( \xi^+)&=& d^+ \geq d_k(\xi).
\end{eqnarray}
This shows (by the
argument at the beginning of the proof and the discussion at the end
of the previous section)
%
\begin{equation} \label{loew}
M( \xi^+,\theta) \geq M(\xi,\theta).
\end{equation}
Moreover, it follows from Chapter II, Section 6 of \citet{karstu1966}
that the point $\vec{d_k}( \xi^+)$ is a boundary point of the ``moment
space''
\[
\mathcal{M}_k = \{ \vec{d_k} (\eta) \mid\eta\mbox{ design on }
[A,B] \}.
\]
Consequently, we obtain from Theorem 2.1 in \citet{karstu1966} that the
design $\xi^+$ is based on at most $\frac{k+2}{2}$ support points,
which proves the first part of the statement.

We now consider the cases (1a) and (1b). The vector $\vec{d}_{k-1}(\xi
)$ is either a~boundary point or an interior point of the $(k-1)$th
moment space $\mathcal{M}_{k-1}$. The first case is characterized by
an index satisfying $I(\xi) < k/2 $ and there exists a unique measure
$\tilde\xi$ with ``moments'' up to the order $k$ specified by $\vec
{d}_{k-1}(\xi)$. To prove this statement regarding uniqueness suppose that
$ I(\xi)< {k\over2}$\vspace*{-1pt} and that there exists a further design, say
$\tilde\xi$, with this property.
{A simple counting argument shows that the total number of distinct
points, say $x_1, \ldots,x_t$ among the support points of both
representations is at most $k$. If it would be less than $k$ we could
take additional support points with corresponding vanishing weights and thus
without less of generality, we can assume that the number of distinct
points is equal to $k$.
Therefore, there would exist $k$ different points
\[
A \leq x_0< x_1<\cdots< x_{k-1}\leq B
\]
such that
\[
\Psi\mu=0,
\]
where the matrix $\Psi$ is given by
\[
\Psi=
\pmatrix{
\Psi_0(x_0) &\Psi_0(x_1) &\ldots&\Psi_0(x_{k-1}) \cr
\Psi_1(x_0) &\Psi_1(x_1) &\ldots&\Psi_1(x_{k-1}) \cr
\vdots&\vdots& \ddots&\vdots\cr
\Psi_{k-1}(x_0) &\Psi_{k-1}(x_1) &\ldots&\Psi_{k-1}(x_{k-1})
}
\]
and the vector $\mu\not= 0$ has components
\[
\mu_i= \cases{
\omega_i,&\quad $x_i \in\operatorname{supp}\xi,
x_i\notin\operatorname{supp}\tilde\xi$,\vspace*{1pt}\cr
-\tilde\omega_i,&\quad $x_i\notin\operatorname{supp}\xi, x_i\in
\operatorname{supp}\tilde\xi$,\vspace*{1pt}\cr
\omega_i-\tilde\omega_i,&\quad $x_i\in\operatorname{supp}\xi\cap
\operatorname{supp}\tilde\xi$,\vspace*{1pt}\cr
0,&\quad $x_i\notin\operatorname{supp}\xi,
x_i\notin\operatorname{supp}\tilde\xi$
}
\]
(here $\omega_i$ and $\tilde\omega_i$ denote the weights of the
designs $\xi$ and $\tilde\xi$, resp.). Because $\mu\not= 0$ it
follows from here that $\det\Psi=0$
which is impossible by the definition of Chebyshev systems.}
Consequently, a design with moments specified by~\eqref{unique} is uniquely
determined and therefore we
take $\xi^+=\tilde\xi$, which has at most $\frac{k+1}{2}$ support points
[see Theorem 2.1 in \citet{karstu1966}, page 42].

If the index of the design $\xi$ satisfies $I(\xi) \ge k/2 $ it
follows from the discussion in Chapter II, Section 6 in \citet{karstu1966}
that the design $\xi^+$ defined by \eqref{b1} and \eqref{b2} is the
upper principal representation of the vector $\vec{d}_{k-1}( \xi)$, which
means that its index is precisely $\frac{k}{2}$ and its support
includes the point $B$. Note that for this argument we require condition
\eqref{cheb+}.

Consequently, if $k=2m+1$ is odd, the upper principal representation
$\xi^+$ has index $m + {1\over2}$ and precisely $m+1$ support points
including the point $B$. On the other hand, if $k=2m$ is even, $\xi^+$
has $m+1$ support points and the boundary points $A$ and $B$ of the
design interval are support points because the index of the design $\xi
^+$ is $m$.

The proof of part (2) of Theorem \ref{thm3.1} is similar [where the upper
principal representation has to be replaced by the lower principal
representation using condition \eqref{cheb-}] and omitted.
\end{pf}

\begin{remark}\label{rem3.2}
(a) Note that Theorem
2.1 in \citeauthor{karstu1966} [(\citeyear{karstu1966}), Chapter II] refers to moment spaces
corresponding to not necessarily bounded measures and the inclusion of the
constant function in the system under consideration guarantees its
application to a moment space corresponding to probability measures as
required in the proof of Theorem \ref{thm3.1}.
An alternative explanation can be given by the generalized equivalence
theorem as stated in \citet{pukelsheim2006}. It follows from this result
that for an optimal design (with respect to the commonly used criteria)
there exist some constants, say $a_i \in\mathbb{R}$, $i=1,\dots, k$,
such that for all support points of the optimal design the identity
\[
\sum^k_{i=1} a_i \Psi_i (x)= c
\]
is satisfied, where $c$ denotes a constant (e.g., for the
$D$-optimality criterion~$c$ is the number of parameters). Since an optimal
design is admissible, the inclusion of the constant function guarantees
that the index of these designs
is at most $k/2$. Note that this is a sufficient but,
generally speaking, not necessary condition.

(b) Note that it follows from the proof of Theorem \ref{thm3.1} that the
conditions \eqref{b1} and \eqref{b2} imply \eqref{loew}, that is,
the superiority
of the information matrix of the design $\xi^+$ with respect to the
Loewner ordering. In many cases (e.g., polynomial regression models),
the converse direction is also true and in these cases it follows from
the proof of Theorem \ref{thm3.1} that a design $\xi$ with index $I(\xi) <
{k\over2}$ can only be ``improved'' (with respect to the Loewner
ordering of the corresponding information matrices) by itself. In fact
we are
not aware of any case where the converse direction does not hold.

(c) Note also that Theorem \ref{thm3.1} provides a solution to the problem
indicated in the example of the \hyperref[sec1]{Introduction}. In the linear regression model
we have $k=2$, therefore we can use the given design $\xi_1$
(concentrating all observations at
$x=0$) as an ``improvement'' of
$\xi_1$. However, because the index of $\xi_1$ is $1/2 < 1$ the
design $\xi_1$ can only be improved by itself (see the previous
remark). In
particular, there does not exist a design $\xi$ which takes
observations at $x=1$ and improves $\xi_1$ in the sense $M(\xi)\geq
M(\xi_1)$.

(d) It is also worthwhile to mention that{ a design }improving the
given design $\xi$ is not necessarily unique. Consider, for example, again
the linear regression model on the interval $[0,1]$ and the design $\xi
$ which has equal masses at the points $0$ and $3/4$. The information
matrix of $\xi$ is given by
\[
M (\xi)= \pmatrix{
1 & {3 \over8} \vspace*{2pt}\cr
{3\over8 } & {9 \over32}
}.
\]
Now define for any $p \in[{1\over2}, {5 \over8} ] $ a design\vspace*{-1pt}
$\xi_p^+$ with masses $p$ and $1-p$ at the points $0$ and ${ 3\over8(1-p)}$,
respectively. Then it follows that
\[
M (\xi^+_p )= \pmatrix{
1 & \displaystyle \frac{3 }{8}\vspace*{3pt} \cr
\displaystyle \frac{3 }{8} & \displaystyle \frac{ 9 }{64(1-p)}
}
\]
and $M (\xi^+_p ) \ge M (\xi)$ for any $p \in[{1\over2}, {5 \over
8} ] $. Note that the choice $p= {5 \over8} $ gives the upper principal
representation $\xi^+ = \xi^+_{5/8}$ with index $1$ and support
points $0$ and $1$, while for $p \in[ {1 \over2}, {5 \over8})$ we have
index $I(\xi^+_p)=3/2$.
\end{remark}

In the remaining part of this section, we will relate the result of
Theorem \ref{thm3.1} to the recent findings of \citet{yang2010}. Note that---in
contrast to Theorem 1 and 2 of \citet{yang2010}---our Theorem \ref{thm3.1} does
not require the differentiability of the functions $\Psi_j$. Moreover, in
some cases it provides a~better description of the admissible designs.
For a more detailed explanation, we note that a Chebyshev system of
functions $\{ u_0,\dots,u_k \}$
is called an extended Chebyshev system, if and only
if for any $a_0, \dots,a_k \in\mathbb{R}$ with $\sum^k_{i=0} a^2_i
\neq0$ the function
\[
\sum^k_{i=0} a_i u_i (x)
\]
has at most $k$ zeros counted with multiplicities in the interval
$[A,B]$. Note that this definition is equivalent to the definition
given in \citet{karstu1966}.
It is in fact proved in \citeauthor{karstu1966} [(\citeyear{karstu1966}), Section 1.2] for the case
of system $u_i(t) =t^i, i=0,\ldots,n$. And the argument can be applied
for general case. Moreover, by definition, an extended Chebyshev system
is always a Chebyshev system.

A simple way of constructing an extended Chebyshev system is the
following [see \citet{karstu1966}, page 19].
Let $w_0,\dots,w_k$ be functions on the interval $[A,B]$
which are either positive or negative.
We now consider the new functions
%
\begin{eqnarray} \label{uj}
u_0 (x) & = & w_0(x), \nonumber\\
u_1(x) & = & w_0(x) \int^x_A w_1(t_1)\,dt_1 ,\nonumber\\[-8pt]\\[-8pt]
& \vdots& \nonumber\\
u_k(x) & = & w_0(x) \int^x_A w_1(t_1) \int^{t_2}_A w_2(t_2) \cdots
\int^{t_{k-1}}_A w_k(t_k) \,dt_k \cdots \,dt_1. \nonumber
\end{eqnarray}
A direct calculation shows that the Wronskian determinant of the
functions $u_0,\dots,u_k$ is given by
\begin{eqnarray} \label{wron}
W_{x}(u_0, \dots, u_k)
&=& \left|
\matrix{
u_0(x) & u^\prime_0(x) & \cdots& u^{(k)}_0(x) \vspace*{1pt}\cr
u_1(x) & u^\prime_1(x) & \cdots& u^{(k)}_1(x) \cr
\vdots& \vdots& \ddots& \vdots\vspace*{1pt}\cr
u_k(x) & u^\prime_k(x) & \cdots& u^{(k)}_k(x)
}\right|\nonumber\\[-8pt]\\[-8pt]
&=& (w_0(x))^{k+1} (w_1(x))^k \cdots(w_{k-1}(x))^2 w_k(x)\nonumber
\end{eqnarray}
and it is shown in Chapter XI in \citet{karstu1966} that the set $\{
u_0,\dots,u_k \}$ of $k$ times differentiable function is an extended
Chebyshev system if and only if
\[
W_{x}(u_0,\dots,u_k) > 0
\]
for all $x \in[A,B]$. On the other hand, this representation provides
a constructive method for checking if a given system of $k$ times
differentiable functions $\{ u_0,\dots,u_k \}$ is a Chebyshev system
on the interval $[A,B]$. To be precise, define $w_0(x)=u_0(x)$ and
recursively differential operators
%
\begin{eqnarray}
\label{u1} D_jf &=& \frac{d}{dx} \biggl( \frac{f}{w_j} \biggr) ;\qquad  j=0,\dots,k, \\
\label{u2} w_{j+1} &=& (D_j D_{j-1} \cdots D_0) u_{j+1} ;\qquad  j=0,1,\dots
,k-1.
\end{eqnarray}
Consequently, the set $\{ u_0,\dots,u_k \}$ is a
Chebyshev system if the functions~$w_0,\break\dots,w_k$ calculated by \eqref
{u1} and \eqref{u2} are all positive on the interval
$[A,B]$.

\begin{remark}\label{rem3.3}
\citet{yang2010} constructed a triangle array of
functions $\{ f_{l,t} \mid t=1,\dots,k; t \leq l \leq k \}$ from the
functions $\Psi_1,\dots,\Psi_k$ induced by the nonlinear regression
model \eqref{1.1} using the recursion
\[
f_{l,t}(x)= \cases{
\Psi^\prime_l(x), & \quad $t=1,\dots,k$, \vspace*{1pt}\cr
\displaystyle \biggl( \frac{f_{l,t-1}(x)}{f_{t-1,t-1}(x)} \biggr)^\prime,&\quad $ 2 \leq t \leq k ; t
\leq l \leq k$.
}
\]
It is now easy to see that the functions $w_1,\dots,w_k$ obtained
from \eqref{u1} and \eqref{u2} with $w_0=1$, $u_j=\Psi_j$ $ (j=1,\dots,k)$
are precisely the functions $f_{ll}$ defined by \citet{yang2010}. As a
consequence, we will obtain the main result of \citet{yang2010} as a
special case of our Theorem \ref{thm3.1} (note that our assumptions regarding
the differentiability are slightly weaker than in this reference).
\end{remark}

\begin{thm}\label{thm3.4}
Let $\Psi_1,\dots,\Psi_k$ denote the
$k$ different functions in the information matrix \eqref{3.11} corresponding
to the nonlinear regression model which are not equal to the constant
function. Assume that $ \Psi_j$ is $(j+1)$ times continuously
differentiable, define $w_0=1$ and for $j=0,\dots,k-1$
\[
w_{j+1} = D_j D_{j-1} \cdots D_0 \Psi_{j+1}
\]
and assume that condition \eqref{def} is satisfied. If
\[
F(x) = w_1(x) \cdots w_k(x) \neq0
\]
for all $x \in[A,B]$, then for any given design $\xi$ there exists a
design $\tilde\xi$, such that $I(\tilde\xi) \le{k \over2}$
\[
M(\tilde\xi, \theta) \geq M(\xi,\theta).
\]
If the index of the design $\xi$ satisfies $I(\xi) < {k \over2}$
then $\tilde\xi$ is uniquely determined in the class of all designs
$\eta$ with moments specified by \eqref{unique} and coincides with
the design $\xi$. Otherwise [in the case $I(\xi) \ge{k \over2}$] the
following assertions are valid.
\begin{enumerate}[(1a)]
\item[(1a)] If $k$ is odd and $F(x) < 0$ on the interval $[A,B]$, then
the design $\tilde\xi$ has at most $(k+1)/2$ support points and
$\tilde\xi$ can be chosen such that the point~$A$ is a support point.
\item[(1b)] If $k$ is odd and $F(x) > 0$ on the interval $[A,B]$, then
the design $\tilde\xi$ has at most $(k+1)/2$ support points and
$\tilde\xi$ can be chosen such that the point~$B$ is a support point.
\item[(2a)] If $k$ is even and $F(x) < 0$ on the interval $[A,B]$,
then the design $\tilde\xi$ has at most $k/2$ support points.
\item[(2b)] If $k$ is even and $F(x) > 0$ on the interval $[A,B]$,
then the design $\tilde\xi$ has at most $k/2+1$ support points
and
$\tilde\xi$ can be chosen such that the points~$A$ and $B$ are
support points.
\end{enumerate}
\end{thm}

\begin{pf}
Let us define $\Psi_0 (x) =1$ and note that
\[
F(x)= \frac{W_{x}(\Psi_0, \dots, \Psi_k)}{W_{x}(\Psi_0, \dots, \Psi_{k-1})}.
\]

Thus if $F(x)>0$ then condition (\ref{cheb+}) is fulfilled and if $F(x)<0$,
then condition (\ref{cheb-}) is fulfilled. Now Theorem \ref{thm3.4} is an immediate corollary
of Theorem \ref{thm3.1}.
\end{pf}

\begin{remark}\label{rem3.5}
Note that
if the constant function appears among the different functions $\{ \Psi
_{ij} \mid1 \leq i \leq j \leq p \}$ in the information matrix
\eqref{3.11} it is not counted in Theorem \ref{thm3.4} or Theorem 2 of \citet
{yang2010} (see the proof of Theorems 3 and 5--7 in this reference).
\end{remark}

A number of interesting applications of Theorem \ref{thm3.4} are given
in \citet{yang2010}. Note that in all examples considered there
the functions under consideration generate a special type of Chebyshev
systems, namely extended Chebyshev systems that can be generated by
formulas (\ref{b2}). This follows from Remark \ref{rem3.3} and the discussion before
Theorem \ref{thm3.4}. Note that several other interesting examples for the
case of two parameters are given in \citet{yangstuf2009}. All these
examples are based on Lemma 1 from that paper and the conditions of this
lemma are in fact imply that the system of the three functions
(corresponding to different elements of the information matrix) is an extended
Chebyshev system. Thus, these examples can also be considered as
particular cases of Theorem \ref{thm3.1}.

The main advantage of Theorem \ref{thm3.1} consists in the fact that the de la
Garza phenomenon can be established by proving that the system under
consideration is a Chebyshev system. For this purpose, several methods
are available which differ from the approach presented in
\citet{yang2010} and in the next section we will consider an example
illustrating the usefulness of Theorem \ref{thm3.1}.

\section{An application to rational regression models}\label{sec4}

In this section, we present a class of nonlinear regression models
where Theorem \ref{thm3.4} [or Theorem 2 in \citet{yang2010}] is not directly
applicable, but the de la Garza phenomenon can be established by an
application of Theorem \ref{thm3.1}. For this purpose, we consider rational
regression models of the form
%
\begin{equation}\label{4.1}
\eta(x,\theta) = {P(x,\theta_{(1)}) \over Q(x,\theta_{(2)})},
\end{equation}
where
\begin{eqnarray*}
P\bigl(x,\theta_{(1)}\bigr)&=& \theta_1 + \theta_2x +\cdots+ \theta_lx^{(l-1)}, \\
Q\bigl(x,\theta_{(2)}\bigr)&=& 1 + \theta_{l+1}x +\cdots+ \theta_{s+l}x^s
\end{eqnarray*}
are polynomials of degree $l-1$ and $s$, respectively, with
corresponding parameters
\[
\theta_{(1)}=(\theta_1,\dots,\theta_l)^T,\qquad
\theta_{(2)}=(\theta_{l+1},\dots,\theta_{l+s})^T.
\]
It is shown in \citet{hestusun1996} that the information matrix for
this model can be written in the form
\[
M(\xi, \theta)= B(\theta)C(\xi,\theta)B(\theta),
\]
where $\theta= (\theta_1,\dots,\theta_{l+s})^T$, $B$ denotes an
appropriate matrix [see \citet{hestusun1996}], the matrix $C$ is given
by
\[
C(\xi,\theta)= \int^B_A [1/Q^4(x)]h(x)h(x)^T \,d \xi(x),
\]
$h(x)= (1,x, \dots, x^{p-1})^T$ denotes the vector of monomials with
$p=l+s$ and $Q(x) $ is a polynomial of degree $s$. Therefore, it
follows that the different functions in the information matrix are
given by
\[
\Psi_1(x)=1/Q^4(x), \dots, \Psi_k(x)=x^{k-1}/Q^4(x),
\]
where $k=2p-1$. Define $\Psi_0(x)=1$, then it is well known [see
\citet{karstu1966b}] that under the conditions:
\begin{itemize}[(a)]
\item[(a)]$ Q(x)$ does not vanish in the interval $[A,B]$;
\item[(b)]$ [Q^4(x)]^{(2p-1)}$ does not vanish in the interval $[A,B]$
\end{itemize}
the functions $\Psi_0, \Psi_1, \dots, \Psi_{2p-1}$ generate a
Chebyshev system on the interval $[A,B]$ and Theorem \ref{thm3.1} is applicable here.

However, we will give an alternative proof of this property which
yields---as a by-product---a constructive condition under which the
condition (b) is fulfilled.
Assume that $Q^4(x)> 0$ for all $x \in[A,B]$ and note that
a Chebyshev system remains a Chebyshev system after multiplication of
all functions by a positive function. Thus, in order to apply Theorem
\ref{thm3.1} it is sufficient to prove that the functions
\[
1, x, x^2, \dots, x^{2p-2}, - Q^4(x)
\]
generate a Chebyshev system on the
interval $[A,B]$. The following lemma provides a sufficient condition
for this property.

\begin{lem}\label{lem4.1}
Assume that the polynomial $Q(x)$ has only
real roots which are either all smaller than $A$ or larger than $B$. If $s
>l-1$, then the functions
\[
1, x, x^2, \dots, x^{2p-2}, \epsilon Q^4(x),
\]
generate a Chebyshev system on the interval $[A,B]$,
where $\epsilon=+1$ if the roots are smaller than $A$ and $\epsilon
=-1$ if the roots larger than $B$.
\end{lem}

\begin{pf}
Based on the assumptions about $Q(x)$, the polynomial $Q^4(x)$ can be
written as $c\!\prod^{4s}_{i=1}(x-\alpha_i)$, where $\alpha_i$ are not necessary
distinct. Clearly, $(Q^4(x))^{(k)}=c\!\sum_{A_k}\!\prod_{j\in A_k}(x-\alpha
_j)$, where $A_k$ is the set of all possible subsets~of $\{1,\ldots
,4s\}$ with
$4s-k$ elements. Define $x_{\min}$ and $x_{\max}$ as the smallest and
largest root of $Q(x)$, then all derivatives of $Q^4(x)$ of even order
less than $4s-1$ are positive outside of the interval $[x_{\min
},x_{\max}]$. Define $u_0(x)=1,\break u_1(x)=x, \dots, u_{2p-2}(x)=x^{2p-2},
u_{2p-1}(x)=Q^4(x).$ By formulas \eqref{u1} and \eqref{u2}, we can
easily calculate that $w_0(x)=1, w_j(x)=j, j=1,\ldots,2p-2,\break w_{2p-1}(x)=
[Q^4(x)]^{(2p-1)}.$ Thus, if $s > l-1$ it follows that $w_{2p-1}(x)$ is
ne\-gative for $x < x_{\min} $ and positive for $ A > x > x_{\max}$.
Therefore (note that $[Q^4(x)]^{(2p-1)}$ has no roots in the interval
$[A,B]$), we have $w_{2p-1}(x) > 0$ for all $x \in[A,B]$. Now the
assertion of Lemma \ref{lem4.1} follows from the formula for the Wronskian
determinant in \eqref{wron} and the fact that a positive Wronskian
determinant is sufficient for the Chebyshev property of the functions
$u_0,\ldots,u_{2p-1}$.
\end{pf}

The following result is now an immediate consequence of Lemma \ref{lem4.1} and
Theorem \ref{thm3.1} (note that we do not repeat the statement of uniqueness of
the latter result).

\begin{thm}\label{thm4.2}
Consider the rational regression model
(\ref{4.1}). Assume that $s >l-1$ and that the polynomial $Q(x)$ has only real
roots, which are either all smaller than $A$ or larger than $B$. Then\vspace*{1pt}
for any design $\xi$ there exists a design $\tilde\xi$ with at most $p$
support points, such that $M(\xi, \theta) \le M(\tilde\xi, \theta
)$. Moreover:
\begin{itemize}[(1)]
\item[(1)] if the index of $\xi$ satisfies $I(\xi) \geq p - \frac
{1}{2}$ and all roots of the polynomial $Q$ are smaller than A,
then $\tilde\xi$ can be chosen such that the support of $\tilde\xi$
contains the point $A$,
\item[(2)]if the index of $\xi$ satisfies $I(\xi) \geq p - \frac
{1}{2}$ and all roots of the polynomial $Q$ are larger
than B, then $\tilde\xi$ can be chosen such that the support of
$\tilde\xi$ contains the point $B$.
\end{itemize}
\end{thm}

\begin{remark}\label{rem4.3}
(a) Theorem \ref{thm4.2} is an extension of Theorem 5 in \citet{hestusun1996}
who investigated only locally D-optimal designs.

(b) Note that \citet{yang2010} considered the classical weighted
polynomial regression model where the different functions in the information
matrix are given by $\Psi_j(x)=\lambda(x)x^{j-1}$, $j=1,\dots, 2p-1$,
where $\lambda$ is a positive function on the interior of the design space,
which is called efficiency function [see \citet{detttram2010}]. His
findings can be generalized in the following way. If there exists a function
$g(x)$ such that
%
\begin{equation}\label{garz}
\biggl( \frac{d}{dx} \biggr)^j \biggl(\frac
{d}{dx}(\lambda
(x)x^{j-1})/g(x) \biggr) = c_j,\qquad  g(x)>0, x\in[A,B],
\end{equation}
for some constants $c_j \in\mathbb{R} \setminus\{0\}$, $j=1,\dots
,2p-1$, then one can denote
\[
\hat\Psi_1(x)=\int^x_0 g(t)\,dt,\qquad  \hat\Psi_j = \Psi_{j-1},\qquad  j=1,\dots
, 2p-1,
\]
and obtains a system of functions satisfying the assumptions of
Theorem~\ref{thm3.4}. In particular, in Theorem 9 of \citet{yang2010} for the case
$\lambda(x)= \exp(x^2)$ the function $g(x)=\lambda(x)=\exp(x^2)$ is
appropriate, while the case $\lambda(x)=(1-x)^{\alpha+1} (1+x)^{\beta
+1}, \alpha> -1, \beta> -1$ requires the choice $g(x)=(1-x)^{\alpha}
(1+x)^{\beta}$. Moreover, the differential equation \eqref{garz}
shows that there are many other
efficiency functions for which the de la Garza phenomenon in the
weighted polynomial regression model
occurs. For example, if $\lambda(x) = 1/(1+x)^n$, $A> -1, n > 2p-2$
one could use
\[
g(x)= 1/(1+x)^{(n+2)}
\]
and it follows that for the weighted polynomial regression model with
this efficiency function any optimal design can be based on at most $p$
points. However, for the rational model of the form \eqref{4.1} such a
technique seemingly does not work. The alternative way is to prove that
the functions $1,x,\dots, x^{k}, \lambda(x)^{-1}$ generate a
Chebyshev system and to use the new Theorem \ref{thm3.1} to establish the de la Garza
phenomenon. Such a method has been realized for the rational model
\eqref{4.1} in the proof of Theorem \ref{thm4.2}.
\end{remark}

\section*{Acknowledgments}

The authors would like to thank Martina Stein, who typed parts of this
manuscript with considerable technical
expertise. We are also grateful to M. Yang and W. J. Studden for
helpful discussions on an earlier version of this paper. The authors
would also like to thank two anonymous referees for very constructive
comments on an earlier version of this paper.


%

\printaddresses

\end{document}